\documentclass[a4paper,11pt,twoside,reqno]{amsart}

\usepackage[utf8]{inputenc}
\usepackage[plainpages=false,pdfpagelabels=true]{hyperref}
\usepackage{amssymb,amsthm}
\usepackage[margin=1in]{geometry}

\DeclareMathOperator{\sech}{sech}

\newtheorem{Satz}{Theorem}[section]
\newtheorem{Prop}[Satz]{Proposition}
\newtheorem{Lem}[Satz]{Lemma}
\newtheorem{Thm}[Satz]{Theorem}

\theoremstyle{definition}

\newtheorem{Bem}[Satz]{Remark}

\newcommand{\tr}{\operatorname{Tr}}

\newcommand{\C}{{\mathbb{C}}}

\newcommand{\dv}{\text{ }dV}
\newcommand{\s}{\mathbb{S}}

\newcommand{\N}{\mathbb{N}}
\newcommand{\Z}{\mathbb{Z}}
\newcommand{\R}{\mathbb{R}}

\allowdisplaybreaks[1]

\renewcommand{\epsilon}{\varepsilon}

\numberwithin{equation}{section}

\title{Infinite families of harmonic self-maps of ellipsoids in all dimensions}
\author{Volker Branding}
\address{University of Vienna, Faculty of Mathematics\\
Oskar-Morgenstern-Platz 1, 1090 Vienna, Austria\\}
\email{volker.branding@univie.ac.at}

\author{Anna Siffert}
\address{Universität Münster, Mathematisches Institut\\
Einsteinstr. 62\\
48149 M\" unster\\
Germany}
\email{asiffert@uni-muenster.de}

\date{\today}

\thanks{The first author gratefully acknowledges the support of the Austrian Science Fund (FWF) through the project "Geometric Analysis of Biwave Maps" (DOI: 10.55776/P34853).
The second author gratefully acknowledges the
support of
the Deutsche Forschungsgemeinschaft (DFG, German Research Foundation) - Project-ID 427320536 - SFB 1442.}

\begin{document}

\begin{abstract}
We prove that for given $k\in\N$, $k\geq 3$, $d\in\N$
and each $a\in\R^{*}$ with
\[
a^2<4d (d+k-2)(k-2)^{-2}
\]
the ellipsoid $E_a:=\{x\in\R^k\,\lvert\,a^{-2}x_1^2+x_2^2+\ldots+x_k^2=1\}$
admits infinitely many harmonic self-maps.
\end{abstract}

\maketitle

\section{Introduction}
A smooth map $\varphi:(M,g)\rightarrow(N,h)$ between closed Riemannian manifolds 
is said to be \textit{harmonic} if it is a critical point of the energy functional
\begin{align}
\label{energy-0}
E(\varphi)=\frac{1}{2}\int_M|d\varphi|^2\dv_g
\end{align}
and can be characterized by the vanishing of its \emph{tension field}
\(\tau(\varphi):=\tr_g\nabla^{\varphi^\ast TN}d\varphi\).
Here, \(\nabla^{\varphi^\ast TN}\) represents the connection on the vector bundle \(\varphi^\ast TN\).
The harmonic map equation
\begin{align*}
\tau(\varphi)=0
\end{align*}
represents a semilinear second order elliptic partial differential equation
with many applications in analysis, geometry and theoretical physics.
The study of harmonic maps is thus a vast research topic, in the present manuscript we will deal exclusively with the existence of harmonic maps.
Due to the nonlinear nature of the harmonic map equation the existence of harmonic maps can be very complicated to achieve.

In their seminal paper \cite{MR164306}, Eells and Sampson started to study the question whether every homotopy class of maps between Riemannian manifolds admits a harmonic representative.
They proved that if all sectional curvatures of the target manifold are non-positive then the answer to this question is affirmative.
However, if the target manifold also admits positive sectional curvatures the answer to this question is only known in special cases.
In this case, in addition to existence results, there are also non-existence results such as that of
Eells and Wood \cite{MR420708} who proved that there is no harmonic map of degree one from a two-torus $\mathbb{T}^2$ onto $\s^2$.

There are many ways to approach the existence problem for harmonic maps to targets of positive curvature,
such as spheres, and it turns out that methods from dynamical systems are a rich source for
this particular kind of problem. We will now provide a short history of this particular direction of research.

The first author who made use of this strategy is Smith: In his seminal article \cite{MR391127} he showed that homotopy classes of spheres up to dimension seven can be 
represented by harmonic maps and to prove this result he used a shooting method for
the corresponding ordinary differential equations.

Later Baird \cite{MR716320} and Baird and Ratto \cite{MR1132860} investigated how the result
of Smith is affected by a deformation of the metric.

In 1990 Eells and Ratto used variational methods to study harmonic maps between ellipsoids
\cite{MR1044657}, see also their book on this subject \cite{MR1242555}.

In 1991 Eells and Ferreira \cite{MR1122903} showed that, 
provided the domain has dimension greater than or equal to three, 
there is always a conformal deformation of the domain metric with respect to which
there is a harmonic representative by considering a modified energy functional and
applying the Palais-Smale condition.

In 1995 Bizo\'{n} \cite{MR1369057} constructed infinitely many harmonic self-maps of $\s^3$, this result was generalized in 1997 to $\s^k$ for $3\leq k\leq 6$ by Bizo\'{n} and Chmaj in \cite{MR1436833}. Both of these articles employ a dynamical systems approach and do not build on tools from the calculus of variations.

In Section\,6 of \cite{MR1369057} Bizo\'{n} speculates that the actual reason for the existence of these infinite families of harmonic maps in specific dimensions could be a topological reason:
\lq The dynamical system approach which I used in this paper gives little insight to what is the actual \lq reason of existence\rq\, of solutions, by which I mean some general mechanism (preferably topological) explaining the existence of solutions in a manner largely insensitive to the details of the model.\rq\,
However, our main result, Theorem\,\ref{thm:main}, shows that this is not the case, i.e. the actual \lq reason of existence\rq\, for these infinite families of harmonic self-maps for specific dimensions of the spheres
is not of topological nature.

\smallskip

In this manuscript we focus on harmonic self-maps of Euclidean ellipsoids.
Although it seems that harmonic maps to ellipsoids seem to have a rich structure
there are not many articles in the literature that have taken up this direction of search.

Below we parametrize the \(k\)-dimensional ellipsoid $E_a\subset\R^k$ as follows
\begin{align*}
x_1=a\cos\psi,~~x_j=z_j \sin\psi, \qquad j=2,\ldots,k, \qquad\text{with}\qquad \sum_{j=2}^kz_j^2=1.
\end{align*}
Further, we assume that $E_a$ is endowed with the induced metric
\begin{align*}
ds^2=(a^2\sin^2\psi+\cos^2\psi)\,d\psi^2+\sin^2\psi\,d\Omega_{k-1},
\end{align*}
where \(d\Omega_{k-1}\) represents the volume element of
the sphere of dimension \(k-1\). Below we parametrize $\s^{k-1}$ by the variable $\theta$.
In this article we study self-maps of $E_a$ of the form
\begin{align*}
(\psi,\theta)\mapsto (f(\psi),\chi(\theta)),
\end{align*}
where $f:[0,\pi]\rightarrow\R$ and $\chi:\s^{k-1}\rightarrow\s^{k-1}$ is an eigenmap of $\s^{k-1}$.
It is straightforward to verify and well-known that a map \(u\colon M\to\s^n\) is harmonic if it solves
\begin{align*}
\Delta u+|du|^2u=0.
\end{align*}
A harmonic map \(u\colon\s^{p-1}\to\s^{q-1},p,q\geq 2\) is called \emph{eigenmap} if 
\(|du|^2=\lambda\) where \(\lambda\) is a constant.
It is well-known that \(u\) is a harmonic eigenmap if and only if all of its
components are harmonic polynomials of common degree \(d\), in which case
\(\lambda_d=d(d+p-2)\). The integer \(\lambda_d\) is called \emph{eigenvalue} of
the eigenmap \(u\). For more details on eigenmaps we refer to the seminal article \cite{MR391127}  of Smith and the book \cite{MR1242555} of Eells and Ratto.

\smallskip

Our first main result is the following existence result.
\begin{Thm}
\label{thm:main}
Let $k\in\N$, $k\geq 3$, and $d\in\N$ be given.   
For each $a\in\R^{*}$ with
\begin{align}
\label{eq:smallness-a}
a^2< 4d (d+k-2)(k-2)^{-2}
\end{align}
the ellipsoid $E_a\subset\mathbb{R}^k$ admits infinitely many harmonic self-maps.
\end{Thm}

On the other hand, our second main result shows that the condition \eqref{eq:smallness-a}
is crucial for our existence result.

\begin{Thm}
\label{thm:non-existence}
Let $k\in\N$, $k\geq 3$, and $d\in\N$ be given.   
If $a\in\R^{*}$ satisfies
\begin{align*}
a^2\geq 4d (d+k-2)(k-2)^{-2}
\end{align*}
the harmonic self-maps of the ellipsoid $E_a\subset\mathbb{R}^k$ of the form as in Theorem \ref{thm:main} do not exist.
\end{Thm}

In the following remarks we put our main result Theorem \ref{thm:main} into its mathematical context.

\begin{Bem}
\begin{enumerate}
\item 
Note that for \(a=1\) the ellipsoid \(E_a\) degenerates to the sphere \(\s^k\).
If in addition we assume $\chi=\mbox{id}$, i.e. $d=1$, we have exactly the case studied by 
Bizo\'{n} and Chmaj in \cite{MR1436833}.
It is straightforward to check that in this case the smallness condition \eqref{eq:smallness-a}
is only satisfied if \(k\leq 6\) which is precisely what one would expect
from the existence result of Bizo\'{n} and Chmaj from \cite{MR1436833}.
\item One can argue that there is a geometric reason why the ellipsoid admits harmonic self-maps
in all dimensions \(k\) while on the sphere we need to impose the condition \(3\leq k\leq 6\).
Both the ellipsoid and the sphere have positive curvature, but for \(a\) being small enough,
one of the sectional curvatures of the ellipsoid also becomes as small as needed. Usually, positive curvature on the target manifold is an obstruction against the existence of harmonic maps. However, one can argue that part of our ellipsoids
are close to flat space such that their curvature no longer prohibits the existence of harmonic maps.

\item The advantage of including the eigenmap is the following: If the eigenmap
is the identity, i.e. $\chi=\mbox{id}$ and $d=1$, the smallness condition \eqref{eq:smallness-a}
reads
\begin{align*}
a^2< 4(k-1)(k-2)^{-2}
\end{align*}
forcing \(a^2\) to be small. By also taking into account the eigenmap we can choose \(a^2\) arbitrarily
as long as we choose \(d\) large enough such that \eqref{eq:smallness-a} is satisfied.

\item The first instance where harmonic maps to ellipsoids were considered
in the mathematics literature seems to be in the seminal article of Smith \cite{MR391127}
on harmonic maps between spheres. In Section 9 he points out a number
of interesting remarks showing that harmonic maps to ellipsoids can be very
different compared to the case of a spherical target.

\item Using variational methods and equivariant differential geometry Eells and Ratto \cite{MR1044657}
constructed harmonic maps between ellipsoids and spheres. In Corollary 5.8 they also considered
the case of harmonic self-maps of ellipsoids and they established the existence of a solution
under the condition \eqref{eq:smallness-a}.
Although their method of proof is completely different and relies on variational methods they need to impose the same smallness condition \eqref{eq:smallness-a} that we need to require in Theorem \ref{thm:main}.
However, by our approach, using methods from dynamical systems, we are able to establish the existence of an
infinite family of solutions instead of just a single solution.

\item In \cite{MR733770} Baldes established the existence of a rotationally symmetric harmonic map
from the Euclidean ball \(B\) of dimension \(k\) taking values in the ellipsoid \(E_a\) with prescribed
boundary data on \(\partial B\) assuming the smallness condition \eqref{eq:smallness-a} with \(d=1\).
The result of Baldes was later extended by Fardoun to the case
of \(p\)-harmonic maps \cite{MR1614611} in which case the smallness condition \eqref{eq:smallness-a}
becomes \(a^2<4(k-1)(k-p)^{-2}\). 
\item An interesting approach to the study of harmonic maps between ellipsoids was presented
in the book of Baird \cite[Subsection 9.3]{MR716320} by considering harmonic
maps between spheres and then performing a suitable deformation of the metric.
\end{enumerate}
\end{Bem}

For the sake of completeness we want to mention the articles of
Helein \cite{MR924091} and Hong \cite{MR1935488} where questions of regularity
and stability of harmonic maps to ellipsoids are discussed.

For a general introduction concerning harmonic maps with symmetries
we refer to the by now classic book of Eells and Ratto \cite{MR1242555}.

\smallskip

Let $d\in\N$ be given.
To prove Theorem\,\ref{thm:main} we reduce the construction of harmonic self-maps of $E_a$ to finding solutions $f:[0,\pi]\rightarrow\R$ of the ordinary differential equation
\begin{align}
\label{el}
f''(\psi)\notag=&\frac{d(d+k-2)}{2}\frac{a^2\sin^2\psi+\cos^2\psi}{a^2\sin^2(f(\psi))+\cos^2(f(\psi))}
\frac{\sin(2f(\psi))}{\sin^2\psi}-(k-1)\cot\psi f'(\psi)
 \\
&-\frac{a^2-1}{2}\frac{\sin(2f(\psi))}{a^2\sin^2(f(\psi))+\cos^2(f(\psi))}f'(\psi)^2
+\frac{a^2-1}{2}\frac{\sin 2\psi}{a^2\sin^2\psi+\cos^2\psi}f'(\psi),
\end{align}
which satisfy the boundary conditions
\begin{align*}
\lim_{\psi\rightarrow 0}f(\psi)=0\quad\mbox{and}\quad \lim_{\psi\rightarrow \pi}f(\psi)=\ell\pi
\end{align*}
for some $\ell\in\Z$. 
Solutions to this boundary value problem are then constructed by using methods from dynamical systems
as was done in \cite{MR1369057} and \cite{MR1436833}.
We would like to point out that the existence of oscillating solutions of \eqref{el} is crucial for the
proof of Theorem \ref{thm:main}.

\smallskip

In order to prove Theorem\,\ref{thm:non-existence} we study the linearization of the Euler Lagrange equation around the critical point 
\(f=\frac{\pi}{2}\) and show that this equation no longer
admits oscillating solutions if the smallness condition \eqref{eq:smallness-a} is violated.

\medskip

The stability of a given harmonic is an important property that characterizes its qualitative behavior.
The intuition is that if a harmonic map is stable, then there does not exist
a second harmonic map \lq nearby\rq, meaning that the critical points of the energy functional \eqref{energy-0} are isolated.
Here we only study the stability of
the identity map when interpreted
as a harmonic self-map of an ellipsoid as it is the only explicit solution of \eqref{el}
that is known.
More details concerning the stability of
harmonic maps can e.g. be found in the introduction of our recent article \cite{MR4477489}.

\smallskip

We prove that  the identity map, considered as a harmonic map of the ellipsoid $E_a$, is unstable.
This is consistent with the expectation that positive curvature is an obstruction to stability.

\begin{Thm}
\label{thm:stability-identity}
For any $k\in\N$, $k\geq 3$, and any $a\in\R^{*}$
the first eigenvalue of the Jacobi operator associated with harmonic self-maps
of ellipsoids evaluated at the identity map is given by
\begin{align}
\label{eq:first-ev}
\lambda_1=a^2(2-k).
\end{align}
Hence, the identity map, when considered as a harmonic self-map
of ellipsoids, is unstable.
\end{Thm}

\begin{Bem}
\begin{enumerate}
\item As one should expect, for \(a=1\), the first eigenvalue \eqref{eq:first-ev}
of the Jacobi operator yields the corresponding result for the sphere, see for example \cite[Section 5]{MR1436833}.
 \item It might be impossible to explicitly calculate the complete spectrum 
of the Jacobi operator associated with the identity map considered as
a harmonic self-map of an ellipsoid. In our recent article \cite{MR4477489} 
we were able to determine such spectra in the case of harmonic self-maps of cohomogeneity one manifolds.
However, we were not able to apply these techniques in the case of the ellipsoid. 
\end{enumerate}

\end{Bem}

\medskip

\textbf{Organisation:}
In Section\,\ref{sec-1} we provide preliminaries.
The proof of Theorem\,\ref{thm:main} is contained in Section\,\ref{sec-2}. We prove Theorem\,\ref{thm:non-existence} in Section\,\ref{sec-3}.
Finally, Theorem\,\ref{thm:stability-identity} is shown in Section\,\ref{sec-4}.

\section{Preliminaries}
\label{sec-1}
In this section we derive the Euler-Lagrange equation \eqref{el}.
Further, we perform a change of variables and provide a Lyapunov function for the resulting differential equation.

\smallskip

Recall that we parametrize the \(k\)-dimensional ellipsoid $E_a\subset\R^k$ as follows
\begin{align*}
x_1=a\cos\psi,~~x_j=z_j \sin\psi, \qquad j=2,\ldots,k, \qquad\text{with}\qquad \sum_{j=2}^kz_j^2=1
\end{align*}
and that we assume that $E_a$ is endowed with the induced metric
\begin{align*}
ds^2=(a^2\sin^2\psi+\cos^2\psi)\,d\psi^2+\sin^2\psi\,d\Omega_{k-1},
\end{align*}
where \(d\Omega_{k-1}\) represents the volume element of
the sphere of dimension \(k-1\).
The energy of the self-maps of $E_a$ of the form
\begin{align*}
(\psi,\theta)\mapsto (f(\psi),\chi(\theta)),
\end{align*}
where $f:[0,\pi]\rightarrow\R$ and $\chi:\s^{k-1}\rightarrow\s^{k-1}$ is an eigenmap of $\s^{k-1}$ with eigenvalue $\lambda=d(d+k-2)$, $d\in\N$,
is given by

\begin{align}
\label{energy}
E(f)=C\int_{0}^{\pi}\bigg(f'^2\frac{a^2\sin^2 f+\cos^2f}{a^2\sin^2\psi+\cos^2\psi}
+d(d+k-2)\frac{\sin^2f}{\sin^2\psi}\bigg)(a^2\sin^2\psi+\cos^2\psi)^\frac{1}{2}
\sin^{k-1}\psi d\psi
\end{align}
for some $C\in\mathbb{R}_+$.
Here and below $'$ denotes derivatives with respect to the variable $\psi$.
The argument of $f$ and its derivatives is always $\psi$ and is therefore henceforth omitted.

\smallskip

The following proposition provides the Euler-Lagrange equation of \eqref{el}. 

\begin{Prop}
The critical points of the energy functional \eqref{energy} are characterized by
\begin{align*}
f''=&\frac{d(d+k-2)}{2}\frac{a^2\sin^2\psi+\cos^2\psi}{a^2\sin^2 f+\cos^2f}
\frac{\sin 2f}{\sin^2\psi}
-\frac{a^2-1}{2}\frac{\sin 2f}{a^2\sin^2f+\cos^2f}f'^2 \\
\nonumber &
-(k-1)\cot\psi f'+\frac{a^2-1}{2}\frac{\sin 2\psi}{a^2\sin^2\psi+\cos^2\psi}f'.
\end{align*}
\end{Prop}

\begin{proof}
We consider a one-parameter variation $f_s$ of \(f\) with fixed endpoints satisfying \(\frac{df}{ds}\big|_{s=0}=\eta\).
A straightforward calculation yields
\begin{align*}
\frac{dE(f_s)}{ds}\big|_{s=0}=C\int_0^{\pi}
\big(2f'\eta'&\frac{a^2\sin^2 f+\cos^2f}{a^2\sin^2\psi+\cos^2\psi}
+\eta d(d+k-2)\frac{\sin 2f}{\sin^2\psi}\\
&+\eta (a^2-1)\frac{\sin 2f}{a^2\sin^2\psi+\cos^2\psi}f'^{2}
\big)(a^2\sin^2\psi+\cos^2\psi)^\frac{1}{2}
\sin^{k-1}\psi d\psi.
\end{align*}
Using integration by parts we deduce
\begin{align*}
\frac{dE(f_s)}{ds}\big|_{s=0}=&C\int_0^{\pi}
\eta\big(-2f''\frac{a^2\sin^2 f+\cos^2f}{a^2\sin^2\psi+\cos^2\psi}
+d(d+k-2)\frac{\sin 2f}{\sin^2\psi}\big) \\
&\times (a^2\sin^2\psi+\cos^2\psi)^\frac{1}{2}
\sin^{k-1}\psi d\psi \\
&-C\int_0^{\pi}2\eta f'
\big(
\frac{a^2\sin^2f+\cos^2f}{(a^2\sin^2\psi+\cos^2\psi)^\frac{1}{2}}\sin^{k-1}\psi
\big)'d\psi.
\end{align*}

A direct calculation shows the following identity
\begin{align*}
\bigg(\frac{a^2\sin^2f+\cos^2f}{(a^2\sin^2\psi+\cos^2\psi)^\frac{1}{2}}&\sin^{k-1}\psi\bigg)' \\
=&\bigg(\frac{\sin 2f(a^2-1)}{a^2\sin^2\psi+\cos^2\psi}f'
-\frac{1}{2}\frac{a^2\sin^2f+\cos^2f}{(a^2\sin^2\psi+\cos^2\psi)^2}\sin 2\psi(a^2-1) \\
&+(k-1)\cot\psi\frac{a^2\sin^2f+\cos^2f}{a^2\sin^2\psi+\cos^2\psi}
\bigg)(a^2\sin^2\psi+\cos^2\psi)^\frac{1}{2}\sin^{k-1}\psi.
\end{align*}
The claim follows from combining the preceding equations.
\end{proof}

\begin{Bem}
For \(a=1\), \eqref{el} reduces to the equation studied by Bizo\'{n} and Chmaj in \cite{MR1436833}, see equation (2.3) therein.
The construction of solutions of \eqref{el} in this manuscript is inspired by the method used in \cite{MR1436833}.
However, we note that for $a\neq 1$ the Euler-Lagrange equation \eqref{el} has a more complicated structure compared to the case $a=1$, e.g.
a term which is quadratic in $f^{'}$ appears in equation \eqref{el}. 
Therefore, the proofs become substantially more demanding.
\end{Bem}

Now, we make the change of coordinates
\begin{align*}
x=\log(\tan\frac{\psi}{2}),\qquad h=f-\frac{\pi}{2}.
\end{align*}
The argument of $h$ and its derivatives is always $x$ and is therefore henceforth omitted.

In terms of these variables the energy now acquires the form
\begin{align*}
E(h)=C\int_{-\infty}^\infty\big(
h'^2\frac{a^2\cos^2 h+\sin^2h}{a^2\mbox{sech}^2x+\tanh^2x}+d(d+k-2)\cos^2h
\big)
(a^2\mbox{sech}^2x+\tanh^2x)^\frac{1}{2}\mbox{sech}^{k-2}xdx.
\end{align*}

The critical points of \(E(h)\)
 are characterized by the following differential equation
\begin{align}
\label{euler-lagrange-h}
h''\nonumber&+\frac{1}{2}(1-a^2)\frac{\sin(2h)}{a^2\cos^2h+\sin^2h}h'^2
-(1-a^2)\frac{\tanh x}{a^2+\sinh^2x}h'
-(k-2)\tanh xh' \\
&+\frac{d(d+k-2)}{2}\frac{a^2\mbox{sech}^2 x+\tanh^2x}{
a^2\cos^2h+\sin^2h}\sin(2h)=0.
\end{align}

From now on we assume $h$ to be a solution of \eqref{euler-lagrange-h}.
We introduce \(W\colon\R\to\R\) by 
\begin{align}
\label{def-lyapunov}
W(x):=h'^2\frac{a^2\cos^2 h+\sin^2h}{a^2\mbox{sech}^2x+\tanh^2x}+d(d+k-2)\sin^2h.
\end{align}
In the following lemma we show that \eqref{def-lyapunov} is a Lyapunov function for \eqref{euler-lagrange-h}.

\begin{Lem}
\label{lemma-w}
The function \(W\colon\R\to\R\) satisfies
\begin{align}
\label{w-derivative}
W'(x)=2(k-2)\tanh xh'^2\frac{a^2\cos^2 h+\sin^2h}{a^2\operatorname{sech}^2x+\tanh^2x}.
\end{align}
In particular, $W$ is monotonically decreasing on $(-\infty,0]$ and monotonically increasing on $[0,\infty)$.
\end{Lem}

\begin{proof}
A direct calculation shows that
\begin{align*}
W'(x)=&2h''h'\frac{a^2\cos^2 h+\sin^2h}{a^2\operatorname{sech}^2x+\tanh^2x} +(1-a^2)\frac{\sin(2h)}{a^2\operatorname{sech}^2x+\tanh^2x}h'^3\\
&-2(1-a^2)\frac{\tanh x}{a^2+\sinh^2 x}\frac{a^2\cos^2 h+\sin^2h}{a^2\operatorname{sech}^2x+\tanh^2x}h'^2 +d(d+k-2)\sin(2h) h'.
\end{align*}
Now, replacing \(h''\) in the above equation using \eqref{euler-lagrange-h} completes the proof.
\end{proof}

\section{Proof of Theorem \,\ref{thm:main}}
\label{sec-2}
In this section we provide the proof of Theorem\,\ref{thm:main}.

\smallskip

From \eqref{euler-lagrange-h} we get that solutions of \eqref{euler-lagrange-h} with $h'(0)=0$ are even under the reflection \(x\to-x\) and solutions with $h(0)=0$ are odd under the reflection \(x\to-x\). Below we will focus on such odd and even solutions only. Therefore it is sufficient to consider $x\geq 0$.
An odd solution of \eqref{euler-lagrange-h} with $h'(0)=b$ will be called \textit{\(b\)-orbit} and denoted by $h(x,b)$.
An even solution of \eqref{euler-lagrange-h} with $h(0)=d$ will be called \textit{\(d\)-orbit}.
Below we focus on \(b\)-orbits, but all our considerations can easily be adapted to \(d\)-orbits as well.
The Lyapunov function (\ref{def-lyapunov}) associated with a \textit{\(b\)-orbit} will be denoted by $W(x,b)$.

\smallskip

The following lemma is an immediate consequence of Lemma\,\ref{lemma-w}. It shows that $W(\cdot,b)$ is small on intervals of the form $[0,T]$, $T\in\R^{+}$, provided that $b$ is small enough.

\begin{Lem}
\label{lemma-close}
Given any $T>0$ and $\eta >0$, there exists an $\epsilon(\eta,T)$ such that if $b<\epsilon$ then $W(x,b)<\eta$ for $x\leq T$.
\end{Lem}
\begin{proof}
From Lemma\,\ref{lemma-w} we have
\begin{align*}
W'(x,b)\leq2(k-2)W(x,b).
\end{align*}
Integrating once hence yields
\begin{align*}
W(x,b)\leq W(0,b)\exp(2(k-2)x)=b^2\exp(2(k-2)x).
\end{align*}
For $\epsilon= \exp(-(k-2)T)\sqrt{\eta}$ the claim hence follows.
\end{proof}

The intuition for Lemma\,\ref{lemma-close} is as follows:
\begin{itemize}
\item
For $a^2\leq 1$ and $\eta$ small, the estimate 
\begin{align*}
 a^2h'^2(x,b)+d(d+k-2)\sin^2h(x,b)\leq W(x,b) < \eta
\end{align*}
and Lemma\,\ref{lemma-close} imply that $h(x,b)$ and $h'(x,b)$ are small for $0\leq x\leq T$.
In other words, the graph of $h$ stays close to $0$ for $0\leq x\leq T$.\\
\item For $a^2\geq 1$ and $\eta$ small, the estimate 
\begin{equation*}
 h'^2(x,b)+d(d+k-2)\sin^2h(x,b)\leq W(x,b) < \eta
\end{equation*}
and Lemma\,\ref{lemma-close} imply that $h(x,b)$ and $h'(x,b)$ are small for $0\leq x\leq T$.
In other words, the graph of $h$ stays close to $0$ for $0\leq x\leq T$.
\end{itemize}

\smallskip

Next, we introduce the so-called rotation number of a $b$-orbit. We will eventually prove that for each $\ell\in\N$ there exists a solution of the boundary value problem \eqref{el} with rotation number $\ell$.\\
For any \(b\)-orbit we set \(\theta(x,b)\) to be
\begin{align*}
\theta(0,b)=\frac{\pi}{2},\qquad \theta(x,b)=\arctan\big(\frac{h'(x,b)}{h(x,b)}\big)
\end{align*}
for any \(x>0\). The \textit{rotation number} \(\Omega(b)\) of the $b$-orbit is defined by 
\begin{align*}
\Omega(b)=-\frac{1}{\pi}\big(\theta(x_e(b),b)-\theta(0,b)\big).
\end{align*}
Here, $x_e(b)$ denotes the smallest $x>0$ at which the $b$-orbit exits the set
\begin{align*}
\Gamma:=\{(h,h',x)\mid h<\frac{\pi}{2}, x>0, (h,h')\neq (0,0)\}.
\end{align*}
If the $b$-orbit does not exit $\Gamma$, we set $x_e(b)=\infty$.

\smallskip

In the next lemma we show that $\theta'(\cdot,b)$ is uniformly bounded by a negative constant on an interval of the form $(x_0,T)$, $x_0>0$, provided that $b$ is small enough.

\begin{Lem}
\label{lemma-theta}
Let $k\in\N$ be given. Then, there exist $a\in(0,1)$, $x_0>0$ and $c>0$ such that for any $T>x_0$ there exists an $\epsilon>0$,
such that
$\theta'(x,b)\leq -c$ for all $x\in(x_0,T)$ and any $b\in(0,\epsilon)$.
\end{Lem}
\begin{proof}
Let $k\in\N$ be fixed.
Throughout the proof we repeatedly use the notations $\theta:=\theta(x):=\theta(x,b)$ and $h(x):=h(x,b)$; we omit the arguments when they worsen the readability and do not contribute to the clarity of the proof.
A straightforward calculation using \eqref{euler-lagrange-h} yields
\begin{align*}
\theta'(x)=&-\sin^2\theta+(1-a^2)\frac{\tanh x}{a^2+\sinh^2x}\frac{\sin(2\theta)}{2}
+(k-2)\tanh x \frac{\sin(2\theta)}{2}\\
&-d(d+k-2)\frac{a^2\text{sech}^2x+\tanh^2x}{a^2\cos^2h(x)+\sin^2h(x)}\frac{\sin(2h(x))}{2h(x)}\cos^2\theta\\
&-(1-a^2)\frac{1}{a^2\cos^2h(x)+\sin^2h(x)}\frac{\sin(2h(x))}{2h(x)}h'^2(x)\cos^2\theta.
\end{align*}
We use the identity
\begin{align*}
-\sin^2\theta=-\frac{1+d(d+k-2)}{2}+d(d+k-2)\cos^2\theta+\frac{1-d(d+k-2)}{2}\cos(2\theta)
\end{align*}
and split
\begin{align*}
-&\frac{a^2\text{sech}^2x+\tanh^2x}{a^2\cos^2h(x)+\sin^2h(x)}\frac{\sin(2h(x))}{2h(x)}\\
=&
(a^2\text{sech}^2x+\tanh^2x)\big(\frac{1}{a^2}-\frac{1}{a^2\cos^2h(x)+\sin^2h(x)}\frac{\sin(2h(x))}{2h(x)}\big)-(\text{sech}^2x+\tanh^2xa^{-2}).
\end{align*}
Furthermore, we use
\begin{align*}
\text{sech}^2x+\frac{\tanh^2x}{a^2}=\text{sech}^2x+\tanh^2x(a^{-2}+1-1)=1+\tanh^2x(a^{-2}-1).
\end{align*}
Thus, we obtain
\begin{align}
\label{theta-abl}
\theta'(x)=\notag&-\frac{1}{2}(1+d(d+k-2)+d(d+k-2)\tanh^2x(a^{-2}-1))\\\notag&+\frac{1}{2}(k-2+(1-a^2)\frac{\tanh x}{a^2+\sinh^2x})\lvert\sin(2\theta)\rvert\\&-\frac{1}{2}((k-2+d)d-1+d(d+k-2)\tanh^2x(a^{-2}-1))\cos(2\theta)+\delta_a(x),
\end{align}
where
\begin{align*}
\delta_a(x):=
&-(1-a^2)\frac{1}{a^2\cos^2h(x)+\sin^2h(x)}\frac{\sin(2h(x))}{2h(x)}h^2(x)\sin^2\theta\\
&+d(d+k-2)(a^2\text{sech}^2x+\tanh^2x)\big(\frac{1}{a^2}-\frac{1}{a^2\cos^2h(x)+\sin^2h(x)}\frac{\sin(2h(x))}{2h(x)}\big)\cos^2\theta \\
&+\frac{k-2}{2}(\tanh x\sin(2\theta)-\lvert \sin(2\theta)\rvert)\\
&+\frac{1-a^2}{2}(\frac{\tanh x}{a^2+\sinh^2x}\sin(2\theta)-\frac{\tanh x}{a^2+\sinh^2x}\lvert \sin(2\theta)\rvert).
\end{align*}

As a next step we introduce $A,B:\R\rightarrow\R$ by
\begin{align*}
&A(x)=\frac{1}{2}(k-2+(1-a^2)\frac{\tanh x}{a^2+\sinh^2x}),\\
&B(x)=-\frac{1}{2}((k-2+d)d-1+d(d+k-2)\tanh^2x(a^{-2}-1)),
\end{align*}
and consider the function
\begin{align*}
g(\theta):=A(x)|\sin 2\theta|+B(x)\cos 2\theta.
\end{align*}
Note that by assumption we have \(A(x)>0\) for $x\geq 0$.
The extremal points of \(g(\theta)\) can be characterized by
\begin{align}
\label{extremal-p}
-B(x)+A(x)\frac{\cos 2\theta}{|\sin 2\theta|}=0.
\end{align}
One way of achieving this is by regularizing the absolute value 
(such that it becomes a differentiable function)
in the definition of \(g(\theta)\).
Since $B(x)<0$, we find that \eqref{extremal-p} vanishes for \(\theta^{*}=-\frac{1}{2}\operatorname{arccot}(\frac{B(x)}{A(x)})\).
It is straightforward to check that
\begin{align*}
g(\theta^{*})
=\sqrt{A(x)^2+B(x)^2}.
\end{align*}
Then, for $x\in\R$ we have
\begin{align*}
\theta'(x)\leq -\frac{1}{2}(1+d(d+k-2)+d(d+k-2)\tanh^2x(a^{-2}-1))+\sqrt{A(x)^2+B(x)^2}+\delta_a(x).
\end{align*}
Note that 
\begin{align*}
\lim_{x\rightarrow\infty}&(-\frac{1}{2}(1+d(d+k-2)+d(d+k-2)\tanh^2x(a^{-2}-1))+\sqrt{A(x)^2+B(x)^2})\\
&= -\frac{1}{2}(1+d(d+k-2)a^{-2})+\frac{1}{2a^2}\sqrt{(k-2)^2a^4+(a^2-d(d+k-2))^2}.
\end{align*}
Since $k\in\N$, $k\geq 3$, and $d\in\N$ are given, there
hence exist $a\in(0,1)$, $x_0>0$ and $c>0$ such that 
\begin{align*}
\theta'(x,b)\leq -2c+\delta_a(x)
\end{align*}
for $x\geq x_0$ and any $b\in\R_+$. 
By Lemma\,\ref{lemma-close}, for any $T>x_0$ there exists an $\epsilon>0$ such that $\delta_a(x)<c$ for $x\leq T$ and $b\in(0,\epsilon)$.
This establishes the claim.
\end{proof}

\begin{Bem}
Note that the smallness condition on \(a\), which we applied in the previous proof,
\begin{align*}
-\frac{1}{2}(1+d(d+k-2)a^{-2})+\frac{1}{2a^2}\sqrt{(k-2)^2a^4+(a^2-d(d+k-2))^2}<0
\end{align*}
is equivalent to
\begin{align*}
a^2\leq 4d (d+k-2)(k-2)^{-2}
\end{align*}
which is precisely the condition \eqref{eq:smallness-a} we assumed in Theorem \ref{thm:main}. 
\end{Bem}

From the two preceding lemmas we obtain that the rotation number $\Omega(b)$ becomes arbitrary large provided that $b$ is small enough.

\begin{Prop}
\label{prop-omega}
For any given $N>0$ there exists an $\epsilon>0$ such that for $b\in(0,\epsilon)$ we have $\Omega(b)>N$.
\end{Prop}
\begin{proof}
Let $k\in\N$ be given.
By Lemma\,\ref{lemma-theta} there exist 
$a\in(0,1)$, $x_0>0$ and $c>0$ such that for any $T>x_0$ there exists an $\epsilon_1>0$,
such that
$\theta'(x,b)\leq -c$ for all $x\in(x_0,T)$ and any $b\in(0,\epsilon_1)$.
Hence, for any $T>x_0$ we have
\begin{align*}
\theta(T,b)-\theta(0,b)&=\theta(T,b)-\theta(x_0,b)+\theta(x_0,b)-\theta(0,b)\\&=\int_{x_0}^T\theta'(s,b)ds+\theta(x_0,b)-\theta(0,b) \\
&\leq -c(T-x_0)+\theta(x_0,b)-\theta(0,b).
\end{align*} 
By Lemma\,\ref{lemma-close}, for any given $T>x_0$, there exists an $\epsilon_2>0$ such that $x_e(b)>T$ for $b\in(0,\epsilon_2)$. 
Below let $\epsilon=\mbox{min}(\epsilon_1,\epsilon_2)$ and $b\in(0,\epsilon)$.
We thus get
\begin{align*}
\Omega(b)>-\frac{1}{\pi}[\theta(T,b)-\theta(0,b)]&\geq 
\frac{1}{\pi}[c(T-x_0)-\theta(x_0,b)+\theta(0,b)].
\end{align*}
From \eqref{theta-abl} we have $\theta'(x,b)\leq c_1$ for $x\leq T$, where $c_1$ is some positive constant depending on $k$ and $a$ but not on $b$.
Consequently, we have $\theta(x_0,b)-\theta(0,b)\leq c_1x_0$.
Now, let $T=x_0+\frac{N\pi+c_1x_0}{c}$. This provides the claim.
\end{proof}

The following lemma states that a $b$-orbit which stays in $\Gamma$ and has only finitely many zeros, satisfies the boundary condition $\lim_{x\rightarrow\infty}h(x,b)=\pm\frac{\pi}{2}$.
We omit the proof since it follows along the lines of the proof of Lemma\,2.6 in \cite{MR4579395}.

\begin{Lem}
\label{lem:connecting-orbit}
Let $h$ be a $b$-orbit such that $h(x,b)\in\Gamma$ for all $x\geq 0$. Further, assume that $h(x,b)$ has a finite number of zeros. Then, we have $\lim_{x\rightarrow\infty}h(x,b)=\pm\frac{\pi}{2}$ and $\lim_{x\rightarrow\infty}h'(x,b)=0$.
\end{Lem}

We are now ready to prove our main result, Theorem\,\ref{thm:main}.

\begin{proof}[Proof of Theorem\,\ref{thm:main}:]
Let $d\in\N$ and $k\in\N$ with $k\geq 3$ be given.
We set
\begin{align*}
S_1:=\{b\mid \textrm{$b$-orbit exits } \Gamma \textrm{ via } h=\frac{\pi}{2}
\textrm{ with } \Omega(b)\leq\frac{1}{2}\}.
\end{align*}
We first prove that the set $S_1$ is non-empty:
Recall that for a $b$-orbit we have \(h(0)=0\) and \(h'(0)=b\).
In order to keep formulas shorter, in what follows we write $h(\cdot)$ and $W(\cdot)$ instead of $h(\cdot,b)$ and $W(\cdot,b)$, respectively.
A direct calculation gives
\begin{align*}
W(0)=b^2.
\end{align*}
We chose $b$ such that \(W(0)>d(d+k-2)\). Due to the monotonicity of \(W(x)\) on $[0,\infty)$, see Lemma\,\ref{lemma-w},
we then have \(W(x)>d(d+k-2)\) for all $x\geq 0$.
Therefore, there exists an $\epsilon>0$ such that \(h'(x)>\epsilon>0\) for all \(x>0\).
Consequently the $b$-orbit constructed above exits the set \(\Gamma\)
through \(h=\frac{\pi}{2}\) with \(\Omega(b)<\frac{1}{2}\).
We let \(b_1:=\inf S_1\).
By Proposition \ref{prop-omega} we have \(b_1>0\).
The \(b_1\)-orbit cannot exit the set \(\Gamma\) via \(h=\frac{\pi}{2}\)
as this would also hold true for any \lq nearby orbit\rq\, with \(b<b_1\)
which would contradict the definition of \(b_1\).
Therefore, the \(b_1\)-orbit stays in \(\Gamma\) for all \(x>0\)
and Lemma \ref{lem:connecting-orbit} thus yields \(\Omega(b_1)=\frac{1}{2}\).

\smallskip

To construct the second orbit we define the set \(S_2\) by
\begin{align*}
S_2:=\{b\mid \textrm{$b$-orbit exits } \Gamma \textrm{ via } h=\frac{\pi}{2}
\textrm{ with } \Omega(b)\leq \frac{3}{2}\}.
\end{align*}
In order to complete the proof we show that \(S_2\) is non-empty and then proceed inductively.

\smallskip

By definition of $b_1$, for \(b<b_1\) we have $\Omega(b)>\frac{1}{2}$.
Below we show that for $b<b_1$ still sufficiently close to \(b_1\), we have $\Omega(b)\leq \frac{3}{2}$ and thus \(b\in S_2\).
To accomplish this goal let $x_A>0$ be such that $h'(x_A)=0$ and $0<h(x_A)<\frac{\pi}{2}$.
By choosing $b<b_1$ sufficiently close to $b_1$,  $x_A$ becomes as large as we want.
Furthermore, we let $x_B$ the smallest $x>x_A$ such that $h(x)=0$. Note that both $x_A$ and $x_B$ exist by the definition of $b_1$. 
From the Euler-Lagrange equation \eqref{euler-lagrange-h} we have
\begin{align}
\label{mon}
h'(x)<0\qquad\mbox{for all}\qquad x\in(x_A,x_B).
\end{align}
Our goal is to prove that for $b$ appropriately chosen we have
$W(x_B)> k-2$. 
Thus, from \eqref{mon} and the monotonicity of $W(x)$, we obtain $h'(x)<0$ for all $x\geq x_A$, confirming that the set $S_2$ is non-empty.

\smallskip

Observe that we may assume without loss of generality that
\begin{align}
\label{est-hs}
h'^2(x)\leq d(d+k-2)\qquad\mbox{for all}\qquad x\in(x_A,x_B).
\end{align}
Indeed, if there would exist a $x_0\in(x_A,x_B)$ such that $h'^2(x_0)>d(d+k-2)$, then the monotonicity of $W(x)$ on $[0,\infty)$ yields the existence of an $\epsilon>0$ such that $h'(x)<-\epsilon$ for all $x\geq x_0$.
This in turn implies that $S_2$ is non-empty.

\smallskip

Next we estimate $W(x_B)-W(x_A)$ in two different ways. Afterwards, we combine these estimates to obtain the desired result.
From \eqref{w-derivative} we get
\begin{align}
\label{diff}
W(x_B)-W(x_A)=2(k-2)\int_{x_A}^{x_B}\tanh s\, h'^2(s)\frac{a^2\cos^2 h(s)+\sin^2h(s)}{a^2\operatorname{sech}^2s+\tanh^2s}ds.
\end{align}
Note that for $a^2\leq 1$, we have
\begin{align}
\label{est-a}
a^{-2}\geq\frac{a^2\cos^2 h(s)+\sin^2h(s)}{a^2\operatorname{sech}^2s+\tanh^2s}\geq a^2.
\end{align}
On the other hand, for $a^2\geq 1$, we have
\begin{align*}
\frac{a^2\cos^2 h(s)+\sin^2h(s)}{a^2\operatorname{sech}^2s+\tanh^2s}\geq a^{-2}.
\end{align*}
Below we restrict ourselves to the case $a^2\leq 1$, the case $a^2\geq 1$ can be treated analogously.
From \eqref{diff} and \eqref{est-a} we get
\begin{align}
\label{diff-2}
W(x_B)\notag&\geq W(x_A)+2(k-2)\tanh(x_A)a^2\int_{x_A}^{x_B}h'^2(s)ds\\
&\geq d(d+k-2)\sin^2(h(x_A))+2(k-2)\tanh(x_A)a^2\int_{x_A}^{x_B}h'^2(s)ds.
\end{align}

From (\ref{def-lyapunov}) and the monotonicity of $W(x)$ on $[0,\infty)$ we get
\begin{align*}
W(x)-W(x_A)=&h'^2(x)\frac{a^2\cos^2 h+\sin^2h}{a^2\mbox{sech}^2x+\tanh^2x}\\
&+d(d+k-2)\sin^2h-d(d+k-2)\sin^2(h(x_A))\geq 0.
\end{align*}
Let $x_1\in(x_A,x_B)$ be such that $h(x_1)=\frac{\pi}{4}$. Then, for $x\in [x_1,x_B]$, the preceding inequality yields
\begin{align*}
h'^2(x)\frac{a^2\cos^2 h+\sin^2h}{a^2\mbox{sech}^2x+\tanh^2x}&\geq d(d+k-2)\sin^2(h(x_A))-d(d+k-2)\sin^2h(x)\\&\geq d(d+k-2)\big(\sin^2(h(x_A))-\frac{1}{2}\big).
\end{align*}
Consequently, by \eqref{est-a} we get
\begin{align*}
h'^2(x)\geq a^2d(d+k-2)\big(\sin^2(h(x_A))-\frac{1}{2}\big)
\end{align*}
for $x\in [x_1,x_B]$.
Plugging this inequality into \eqref{diff-2} yields
\begin{align}
\label{est-wb}
W(x_B)\geq \notag&d(d+k-2)\sin^2(h(x_A))\\&+2a^4(k-2)d(d+k-2)\tanh(x_A)(x_B-x_1)\big(\sin^2(h(x_A))-\frac{1}{2}\big).
\end{align}
Since $h(x_1)=\frac{\pi}{4}$, $h(x_B)=0$, we get 
\begin{align*}
x_B-x_1\geq \frac{\pi}{4\sqrt{d(d+k-2)}}
\end{align*}
from \eqref{est-hs}.
Plugging this inequality into \eqref{est-wb} gives
\begin{align*}
W(x_B)\geq &d(d+k-2)\sin^2(h(x_A))\\
&+a^4(k-2)\sqrt{d(d+k-2)}\tanh(x_A)\frac{\pi}{2}\big(\sin^2(h(x_A))-\frac{1}{2}\big).
\end{align*}
For given $a$, we chose $b$ sufficiently close to $b_1$ such that $h(x_A)$ and $\tanh(x_A)$ become close to $\frac{\pi}{2}$ and $1$, respectively.
Consequently, for such $b$ we have $W(x_B)>d(d+k-2)$, which establishes the claim.
This completes the proof of Theorem \ref{thm:main}.
\end{proof}

\section{Proof of Theorem \,\ref{thm:non-existence}}
\label{sec-3}
In this section we provide the proof of Theorem\,\ref{thm:non-existence}.
Crucial for the proof of Theorem\,\ref{thm:main} is the existence of oscillating solutions of \eqref{euler-lagrange-h}.
By studying the linearization of \eqref{euler-lagrange-h} around the critical point 
\(h=0\) we show that this equation no longer
admits oscillating solutions if the smallness condition \eqref{eq:smallness-a} is violated.

\begin{proof}[Proof of Theorem\,\ref{thm:non-existence}:]
In order to show that the Euler-Lagrange equation \eqref{euler-lagrange-h}
does not admit an infinite family of solutions if the smallness condition \eqref{eq:smallness-a}
is not imposed we consider the linearization of \eqref{euler-lagrange-h} around the critical point 
\(h=0\) and show that for large values of \(x\) the linearized equation no longer
admits oscillating solutions.

In order to derive the linearization of \eqref{euler-lagrange-h} we consider
the following one-parameter variation 
\begin{align*}
\frac{dh_s(x)}{ds}\big|_{s=0}=\eta(x).
\end{align*}
In the following we will drop the variable \(x\) and write \(h,\eta\) instead of \(h(x)\)
and \(\eta(x)\) in order not to blow up the notation.

Now, a straightforward calculation shows that the linearization of \eqref{euler-lagrange-h} is given by
\begin{align*}
\eta''&+(1-a^2)\frac{\cos 2h}{a^2\cos^2h+\sin^2h}h'^2\eta
+(1-a^2)\frac{\sin 2h}{a^2\cos^2h+\sin^2h}h'\eta' \\
\nonumber&-\frac{1}{2}(1-a^2)^2\frac{\sin^22h}{(a^2\cos^2h+\sin^2h)^2}h'^2\eta
-(1-a^2)\frac{\tanh x}{a^2+\sinh^2x}\eta'
-(k-2)\tanh x\eta'\\
\nonumber&+d(d+k-2)\frac{a^2\operatorname{sech}^2x+\tanh^2x}{a^2\cos^2h+\sin^2h}\cos(2h)\eta\\
\notag&-(1-a^2)\frac{d(d+k-2)}{2}\frac{a^2\operatorname{sech}^2x+\tanh^2x}{(a^2\cos^2h+\sin^2h)^2}\sin^2(2h)\eta=0.
\end{align*}

Evaluating the above equation at the critical point \(h=0\) yields
\begin{align*}
\eta''-(1-a^2)\frac{\tanh x}{a^2+\sinh^2x}\eta'-(k-2)\tanh x\eta'
+\frac{d(d+k-2)}{a^2}(a^2\operatorname{sech}^2x+\tanh^2x)\eta=0.
\end{align*}
For large values of \(x\) we then find that this equation can be approximated by 
\begin{align*}
\eta''-(k-2)\eta'
+d(d+k-2)a^{-2}\eta=0.
\end{align*}

Making an ansatz \(\eta(x)=Ae^{\alpha x}\) for \(A\in\R,\alpha\in\C\), we obtain
the algebraic equation
\begin{align*}
\alpha^2-(k-2)\alpha+d(d+k-2)a^{-2}=0,
\end{align*}
which has the solutions
\begin{align*}
\alpha=\frac{k-2}{2}\pm\frac{1}{2}\sqrt{(k-2)^2-4d(d+k-2)a^{-2}}.
\end{align*}
This equation has imaginary roots, which is necessary for us in
order to obtain an infinite family of solutions, if and only if
\begin{align*}
a^2<4d(d+k-2)(k-2)^{-2}.
\end{align*}
This is precisely the smallness condition \eqref{eq:smallness-a}.
On the other hand, it is straightforward to see that for
\begin{align*}
a^2\geq 4d(d+k-2)(k-2)^{-2},
\end{align*}
all solutions \(\alpha\) will be real numbers and in this case
we cannot find oscillating solutions of \eqref{euler-lagrange-h}.
\end{proof}

\section{Stability of the identity map}
\label{sec-4}
In this section we prove Theorem \ref{thm:stability-identity}.
As a first step, we provide the Jacobi equation associated with \eqref{euler-lagrange-h}.

\begin{Prop}
\label{prop-stability}
Let \(h\colon\R \to\R\) be a solution of \eqref{euler-lagrange-h}.
Then, the Jacobi equation associated with harmonic self-maps of the ellipsoid 
is given by
\begin{align}
\label{eq:jacobi-general}
\xi''&
+(1-a^2)\frac{\cos 2h}{a^2\cos^2h+\sin^2h}h'^2\xi
+(1-a^2)\frac{\sin 2h}{a^2\cos^2h+\sin^2h}h'\xi'\\
\nonumber &-\frac{1}{2}(1-a^2)^2\frac{\sin^22h}{(a^2\cos^2h+\sin^2h)^2}h'^2\xi
-(1-a^2)\frac{\tanh x}{a^2+\sinh^2x}\xi'
-(k-2)\tanh x\xi' \\
\nonumber&-\frac{d(d+k-2)}{2}(1-a^2)\frac{a^2\sech^2x+\tanh^2x}{(a^2\cos^2h+\sin^2h)^2}\sin^2(2h)\xi\\
\notag&+d(d+k-2)\frac{a^2\sech^2(x)+\tanh^2x}{a^2\cos^2h+\sin^2h}\cos(2h)\xi+\lambda\frac{\xi}{\cosh^2x(a^2\sech^2x+\tanh^2x)}=0,
\end{align}
where \(\xi\colon\R\to\R\) and \(\lambda\in\R\).
\end{Prop}
\begin{proof}
We consider a variation of \(h(x)\) satisfying \(\frac{d}{ds}\big|_{s=0}h_s(x)=\xi(x)\).
Then differentiating \eqref{euler-lagrange-h} with respect to \(s\) yields the result.
Note that the factor multiplying the eigenvalue \(\lambda\) arises
due to the volume element of the ellipsoid.
\end{proof}

In this manuscript we study the stability of the identity map. Recall that for the identity map we have $d=1$.
The following lemma is an immediate consequence of Proposition\,\ref{prop-stability}.

\begin{Lem}
Consider the identity map as a harmonic self-map of the ellipsoid.
Then, the Jacobi equation \eqref{eq:jacobi-general} simplifies to
\begin{align}
\label{jacobi-identity}
\xi''&
-(k-2)\tanh x\,\xi'
+(1-a^2)\frac{\tanh x}{a^2+\sinh^2x}\xi'\\
\nonumber&+\frac{(1-a^2)(a^2-\sinh^2x)}{(a^2+\sinh^2x)^2}\xi 
+\frac{(k-1)(a^2-\sinh^2x)}{a^2+\sinh^2x}\xi
+\frac{\lambda}{\cosh^2x(a^2\sech^2x+\tanh^2x)}\,\xi=0.
\end{align}
\end{Lem}

\begin{proof}
In our setup, the identity map $h_1$ is parametrized by $h_1(x)=-\frac{\pi}{2}+2\arctan(e^x).$
Thus, we have the following identities
\begin{align*}
h_1'=\text{sech}\,x,\quad \cos h_1=\text{sech}\,x,
\quad \sin h_1=\tanh x, \quad \cos 2h_1=2\text{sech}^2x-1,
\quad \sin 2h_1=2\frac{\tanh x}{\cosh x}.
\end{align*}
Plugging these identities into the Jacobi equation, see Proposition\,\ref{prop-stability}, 
we find
\begin{align*}
\xi''-&(k-2)\tanh x\,\xi'
+(1-a^2)\frac{\tanh x}{a^2+\sinh^2x}\xi'
-\frac{1}{2}(1-a^2)^2\frac{4\tanh^2 x}{(a^2+\sinh^2x)^2}\xi \\
&+(1-a^2)\frac{2\sech^2x-1}{a^2+\sinh^2x}\xi 
-\frac{(k-1)}{2}(1-a^2)\frac{4\tanh^2x}{a^2+\sinh^2x}\xi \\
&+(k-1)(2\sech^2x-1)\xi 
+\lambda\frac{\xi}{\cosh^2x(a^2\sech^2x+\tanh^2x)}\,\xi=0.
\end{align*}
Note that we have the following identity
\begin{align*}
-2(1-a^2)\frac{\tanh^2 x}{a^2+\sinh^2x}
+2\sech^2x-1
=&\frac{a^2-\sinh^2x}{a^2+\sinh^2x},
\end{align*}
and applying this twice then completes the proof.
\end{proof}

\begin{proof}[Proof of Theorem\,\ref{thm:stability-identity}]
In order to find the first eigenvalue of the Jacobi operator 
associated with the identity map a
direct calculation shows that
\begin{align*}
\xi(x)=\frac{1}{\sqrt{a^2+\sinh^2x}},\qquad \lambda=a^2(2-k)
\end{align*}
solves \eqref{jacobi-identity} completing the proof.
\end{proof}

\textbf{Data Availability Statement}: Data sharing not applicable to this article as no datasets were generated or
analysed during the current study.

\bibliographystyle{plain}
\bibliography{mybib}
\end{document}